\documentclass[pdflatex,sn-mathphys-num]{sn-jnl}

\usepackage{graphicx}%
\usepackage{multirow}%
\usepackage{amsmath,amssymb,amsfonts}%
\usepackage{amsthm}%
\usepackage{mathrsfs}%
\usepackage[title]{appendix}%
\usepackage{xcolor}%
\usepackage{textcomp}%
\usepackage{manyfoot}%
\usepackage{booktabs}%
\usepackage{algorithm}%
\usepackage{algorithmicx}%
\usepackage{algpseudocode}%
\usepackage{listings}%

\usepackage{bm}
\usepackage{multicol}
\usepackage{caption}
\usepackage{subcaption}
\usepackage[normalem]{ulem}
\usepackage{thm-restate}
\usepackage{xspace}
\usepackage{array}

\usepackage{comment}

\usepackage{enumerate}
\usepackage{enumitem}
\usepackage{mdframed}

\newcommand{\bb}{\mathbb}
\newcommand{\R}{\bb R}
\newcommand{\Z}{{\bb Z}}

\DeclareMathOperator*{\icomp}{icomp}

\def\ve#1{\mathchoice{\mbox{\boldmath$\displaystyle\bf#1$}}
{\mbox{\boldmath$\textstyle\bf#1$}}
{\mbox{\boldmath$\scriptstyle\bf#1$}}
{\mbox{\boldmath$\scriptscriptstyle\bf#1$}}}

\newcommand{\x}{{\ve x}}
\newcommand{\y}{{\ve y}}
\newcommand{\z}{{\ve z}}
\newcommand{\vv}{{\ve v}}
\newcommand{\g}{{\ve g}}
\newcommand{\e}{{\ve e}}
\renewcommand{\u}{{\ve u}}
\renewcommand{\a}{{\ve a}}

\newcommand{\0}{{\ve 0}}

\newcommand{\I}{\mathcal{I}}

\newcommand{\cI}{\mathcal{I}}
\newcommand{\cJ}{\mathcal{J}}
\newcommand{\cO}{\mathcal{O}}

\renewcommand{\epsilon}{\varepsilon}
\newcommand{\ones}{\mathbf{1}}
\newcommand{\cG}{\mathcal{G}}

\let\eps\varepsilon

\def\namedlabel#1#2{\begingroup
    #2%
    \def\@currentlabel{#2}%
    \phantomsection\label{#1}\endgroup
}

\usepackage{physics}

\renewcommand{\ip}[2]{\langle #1, #2 \rangle}

\newcommand{\cQ}{\mathcal{Q}}

\newcommand{\doubleR}{\mathbb{R}}

\newcommand{\remove}[1]{}

\makeatletter
\def\namedlabel#1#2{\begingroup
    #2%
    \def\@currentlabel{#2}%
    \phantomsection\label{#1}\endgroup
}

\theoremstyle{thmstyleone}%
\newtheorem{theorem}{Theorem}%

\theoremstyle{thmstyletwo}%

\theoremstyle{thmstylethree}%
\newtheorem{definition}{Definition}%

\begin{document}

\title[Tight binary first-order oracle lower bounds]{Tight Lower Bounds for Binary First-Order Oracles for Convex Optimization}

\author[1]{\fnm{Amitabh} \sur{Basu}}\email{abasu9@jhu.edu}

\author*[2]{\fnm{Phillip} \sur{Kerger}}\email{kerger@berkeley.edu}

\author[3]{\fnm{Marco} \sur{Molinaro}}\email{mmolinaro@microsoft.com}

\affil[1]{\orgdiv{Department of Applied Mathematics and Statistics}, \orgname{Johns Hopkins University}, \country{USA}}

\affil*[2]{\orgdiv{Department of Industrial Engineering and Operations Research}, \orgname{UC Berkeley}, \country{USA}}

\affil[3]{\orgdiv{Microsoft Research and Computer Science Department}, \orgname{PUC-Rio}, \country{Brazil}}

\abstract{We establish new lower-bounds for the information complexity of mixed-integer convex optimization under two ``bit-wise'' oracles. The first oracle provides bits of first-order information in the standard coordinate model, and the second oracle answers whether the inner product of a specified vector with the gradient of the function at a point or the normal vector of a separating hyperplane for the feasible region is positive or non-positive, thus also providing one bit of first-order information. The new contribution is that under such oracles, the complexity is quadratic in the number of continuous decision variables, which was not known before even for continuous convex optimization. These new lower-bounds are tight (up to a logarithmic term), matched by a natural discretization of standard cutting-plane methods for convex optimization. These reveal that using a standard bit-representation of the first-order information is, in general, the best one can do with respect to the \textit{number of bits of information} needed to solve constrained convex optimization problems.}


\maketitle

\section{Introduction}
We focus on mixed-integer convex optimization problems of the form 
\begin{align}
    \min\{ f(\x, \y): (\x, \y)\in C \cap (\Z^n \times \R^d) \}, \label{def: constrained opti}
\end{align}
where $f: \R^n \times \R^d \rightarrow \R$ is a convex function and $C\subseteq \R^n \times \R^d$ is a closed, convex set. 
The goal is to design algorithms to obtain {\em $\eps$-approximate solutions} for~\eqref{def: constrained opti} for any given error tolerance $\eps > 0$, i.e., return $\z^* \in C \cap (\Z^n \times \R^d)$ such that $f(\z^*) \leq f(\z) + \eps$ for all $\z\in C \cap (\Z^n \times \R^d)$. 

To develop an algorithmic theory for problem~\eqref{def: constrained opti} we must formalize how the instance is presented to the algorithm. This is not an issue if one considers subfamilies of instances with explicit algebraic descriptions (e.g., linear or polynomial optimization). In order to formally handle more general nonlinear instances, the algorithm is given access to the instances via oracles, i.e., the algorithm collects information about the given instance by querying oracles. In this paper, we consider only oracles using first-order information. These are oracles whose answers about the instance $(f,C)$, at a given query point $\z\in \R^n \times \R^d$, are some function of either $f(\z)$, the subgradient of $f$ at $\z$, or of a separating hyperplane to the constraint set $C$ if $\z \not\in C$. For example, the standard first-order oracle returns the entire function value, subgradient, or separating hyperplane. In this paper, we wish to consider more fine-grained oracles that could potentially respond with lesser information about these objects. For example, in any computational setting one always only has access to a bit-representation of subgradients or separating hyperplanes. Given any such oracle access, a fundamental question then is the following: \begin{quote} What is the smallest number of oracle queries one needs to provably report an $\eps$-approximate solution?\end{quote} The smallest number of queries is called the {\em information or oracle complexity} for the given class of instances under that particular oracle. For the standard first-order oracle, this problem has been studied since the mid 70s and in the last few years, a fairly complete understanding has been achieved (which we summarize formally below). Other oracles are much less studied and yet are important to investigate from both a practical and theoretical perspective, connecting with other work on the communication and space complexity of convex optimization.

The feasibility problem of checking whether a closed convex set $C\subseteq \R^n \times \R^d$ contains a point from $\Z^n \times \R^d$ or not, and reporting such a point if one exists is a special case of~\eqref{def: constrained opti} by simply setting $f$ to be the function that takes value $0$ everywhere. Thus, the information complexity of the feasibility problem is a lower bound on the information complexity of~\eqref{def: constrained opti}, for any class of instances and for any oracle. In this paper, we focus on the feasibility problem. In other words, given some class $\I$ of closed, convex sets, and a separation oracle for this class, we wish to understand the smallest number of separation oracle queries needed to either report a feasible point from the given instance, or conclude that the instance is infeasible, i.e., $C = \emptyset$. It is not hard to argue that if $\cI$ is the class of all possible closed, convex sets, then the information complexity is infinite. Indeed, even with $n=0, d=1$ (i.e., detecting feasibility in 1 dimensional continuous convex optimization) for any query point $\hat{y} \in \R$, one could report a separating hyperplane $\{y \in \R: y \leq \hat{y} - 1\}$. After finitely many queries, if the algorithm reports infeasibility, then one can simply consider the instance $\{y : y \leq y_{\min} - 1\}$, where $y_{\min}$ is the smallest $\hat{y}$ queried by the algorithm; on the other hand, if the algorithm reports any point $y^\star$ as a potential feasible point, then one can say that the true instance was a convex set that is contained in $\{y : y \leq y_{\min} - 1\}$ that does not contain $y^\star$. Thus, one has to a priori fix $R > 0$ and consider only instances $C$ such that $C \subseteq [-R, R]^{n + d}$. Without such a bound on the norms of the feasible solutions, the information complexity is clearly infinite. However, this is not enough. By considering all possible singletons in the interval $[0,1]$, a similar argument shows that finitely many queries cannot succeed. Thus, one considers another parameter $\rho > 0$ and focuses on the instances $\cI_{n,d,R,\rho}$ of all closed convex sets $C \subseteq \R^n \times \R^d$ such that $C \subseteq [-R, R]^{n+d}$, and if $C$ is feasible then there exist a ``fat'' set of feasible solutions, i.e., there exists $(\x^\star, \y^\star) \in (\Z^n \times \R^d)$ such that $\{(\x^\star,\y) \in \R^n \times \R^d : \|\y - \y^\star\|_\infty \leq \rho\} \subseteq C$. The information complexity can then be obtained as a function of $n, d, R, \rho$.

Note that for a point $\z\notin C$, the oracle has multiple choices for separating hyperplanes. We will formalize this choice as saying the oracle has a \textit{first-order map} $\cG$ which determines which separating hyperplane is used to answer the query. We give a formal definition here for feasibility problems (see~\cite{Basu2023InformationCO} for the analogous extension for the first-order information of optimization problems).

\begin{definition}\label{def: OUFOI}
A {\em first order map} $\cG$ for $\cI_{n,d,R, \rho}$ is a function $\cG: \cI_{n,d,R,\rho}\times (\R^n \times \R^d) \rightarrow \R^n \times \R^d$ such that $\cG(C, \hat{\z}) = \0$ if $\hat{\z} \in C$; otherwise, $\langle \cG(C, \hat{\z}), \z \rangle < \langle \cG(C, \hat{\z}), \hat{\z} \rangle$ for all $\z \in C$.

An oracle $\cO(\cQ, \cG)$ using first-order information, with \textit{permissible queries} $\cQ$ and \textit{first-order map} $\cG$, takes as input a point $\z \in \R^n \times \R^d$ and a query $q\in \cQ$, which is a function on $\R^n\times \R^d$, and returns for an instance $C \in \cI_{n,d,R, \rho}$ the answer 
$$
\cO(\z, q, C) = q(\cG(C, \z)),
$$
\end{definition}
To summarize: 
\begin{itemize}
    \item Without knowing the instance $C$, based only on the oracle responses thus far, in each step the algorithm defines its query by choosing $q$ from the set of permissible queries $\cQ$ (i.e., what information about the separating hyperplane it is requesting the oracle) and by choosing the query point $\z$. 
    \item The oracle uses $\cG$ to determine \textit{which} separating hyperplane at $\z$ (if any) the oracle uses to answer the query $q$.
\end{itemize}

To define the information complexity of~\eqref{def: constrained opti}, we consider the supremum over all first order maps $\cG$; that is, we consider the information complexity to be the worst case over possible responses that the oracle is permitted to give. With this in mind, we give a formal definition of the information complexity of a class of instances with access to an oracle using first-order information: 

\begin{definition} \label{def: icomp} Given a family of instances $\cI$ and access to an oracle $\cO(\cQ, \cG)$ using first-order information, 
the {\em information complexity} $\icomp(\cI, \cQ)$ of solving the feasibility problem for $\cI$ with access to an oracle of the type in Definition~\ref{def: OUFOI},
is defined as the minimum natural number $k$ 
such that there exist an algorithm that solves the feasibility problem for $\I$ by making at most $k$ queries using the oracle $\cO(\cQ, \cG)$ for any first order map $\cG$. 
\end{definition}

From the lower bound perspective, it thus suffices to describe a way for the oracle to give responses to any algorithm's queries in a way such that after $k$ queries, there are two {\em disjoint} instances, i.e., disjoint sets in $\cI_{n,d,R, \rho}$, that are both consistent with the oracle's answers (i.e., such that there exists a first-order map $\cG$ under which the oracle gives the same answers for these two instances for the queries made by the algorithm). In such a case, the queries made by the algorithm are not able to produce a valid feasible solution (i.e., if it outputs a point in one of these sets, the ``real instance'' could be the other disjoint set). The information complexity is then a function of the class of instances in consideration, and what kinds of queries $\cQ$ on the first-order information can be made.

When $\cQ$ contains only the identity function, $\cO(\cQ, \cG)$ is a standard (exact) first-order separation oracle. However, Definition~\ref{def: OUFOI} captures many other oracles that do not give this full, exact information. In this paper, we establish results for two different such oracles. First, an oracle with queries $q_{ij}$ returning the $i^{th}$ bit of the $j^{th}$ coordinate of the separating hyperplane, and second, an oracle with queries $q_{\vv}$ for $\vv\in \R^d$ returning $1$ if the inner product of $\vv$ and the normal vector $\cG(C, \z)$ for the separating hyperplane is positive, and $0$ otherwise. We will call these the \textit{bit oracle} and the \textit{inner product oracle}, respectively, and denote the set of these queries as $\cQ^{bit}$ and $\cQ^{dir}$. 

This paper focuses on the following question: How many \textit{bits} of first-order information are needed to solve an optimization problem of the form~\eqref{def: constrained opti}? More formally, we study the information complexity of the feasibility version of~\eqref{def: constrained opti} under the bit and inner-product oracles and provide new lower bounds that close the existing gap for these oracles.

\subsection{Our contributions}
We first state the best existing result for binary first-order queries. 
\begin{theorem}\cite[Theorem 9]{Basu2023InformationCO}\label{thm:binary-LB}
For an oracle of the type in Definition~\ref{def: OUFOI} allowing any binary query, i.e., $\cQ$ is the set of all 0/1 functions on $\R^n \times \R^d$ in Definition~\ref{def: OUFOI}, one has   
$$\icomp (\cI_{n,d,R, \rho}, \cQ) = \tilde\Omega \left( 2^n \left(1+\max\left\{d^{\frac{8}{7}}, d \log\left(\frac{R}{\rho}\right)\right\}\right)\right),$$ 

where $\tilde \Omega$ hides polylogarithmic factors in $d$.
\end{theorem}

The above lower bound is obtained using lower bounds on the space complexity of continuous convex optimization in~\cite{marsden2022efficient}, as well as classical lower bounds for the information complexity of continuous convex optimization from Nemirovski and Yudin  \cite{Nemirovski_Yudin79}, as found in \cite[Theorem 9]{Basu2023InformationCO}. In fact, the lower bound of $2^nd \log\left(\frac{R}{\rho}\right)$ can be established for the oracle that returns the full separating hyperplane, as opposed to partial/binary information. These lower bounds are complemented by the following upper bounds.

\begin{theorem}\cite[Theorems 10, 11]{Basu2023InformationCO}
    For both the $\cQ^{bit}$ and $\cQ^{dir}$ oracles, there is an algorithm that solves the feasibility problem for $\cI_{n,d,R, \rho}$ with at most $$O\left(2^nd(n+d)^2\log^2\left(\frac{(n+d)R}{\rho}\right)\right)$$ queries, for any first order map $\cG$. For the continuous problem ($n=0$), this can be improved to $$O\left(d^2\log^2\left(\frac{dR}{\rho}\right)\right)$$
\end{theorem}

These upper bounds are obtained by the natural idea of approximating the separating hyperplane normal to enough bits of accuracy combined with existing cutting-plane algorithms based on exact separation oracles. One immediately notices that, even for the continuous convex optimization case, the upper and lower bounds do not match. The dependence on the number of integer variables is $2^n$ in both bounds; ignoring the logarithmic terms, the gap is in the dependence on the number of continuous variables. The lower bound gives a superlinear $d^{8/7}$ bound, but the naive idea of approximating the separating hyperplane with enough bits of accuracy needs $\Omega(d^2)$ bits in the worst case. This is because even with exact information, these algorithms make $\Omega(d)$ queries in the worst case, and in each query one is using at least $d$ bits of information since there are $d$ coordinates of the separating hyperplane normal vector. The question is: {\em Is this best possible, or is there a more sophisticated query strategy where one does not query all the coordinates in every query and with only subquadratic bits of information one can solve the feasibility (or more generally, the constrained optimization) problem?}

Our main results show that $\Omega(d^2)$ bits are indeed necessary for the bit and inner product oracles.

\begin{theorem} \label{thm:bit-mixed}
Let $d\geq 1$ and $\rho < \frac{R}{2}$. With access to a bit oracle with permissible queries $\cQ^{bit}$, we have
$$\icomp (\cI_{n,d,R, \rho}, \cQ^{bit}) = \Omega \left(2^nd^2 \log\left(\frac{R}{\rho}\right)\right).$$ In fact, this holds even if the query returns the entire coordinate as opposed to just a single bit of a desired coordinate of the normal vector of the separating hyperplane.
\end{theorem}

\begin{theorem} \label{thm:dir-mixed}
Let $d\geq 2$ and $\rho < \frac{R}{6d}$. With access to an inner product oracle with permissible queries $\cQ^{dir}$, we have
$$\icomp (\cI_{n,d,R, \rho}, \cQ^{dir}) = \Omega \left(2^nd^2\left(1 + \frac{1}{\log d}\log\left(\frac{R}{\rho}\right)\right)\right).$$ In fact, this holds even if the query returns the entire inner product of the desired direction and the normal vector of the separating hyperplane.
\end{theorem}

These lower bounds show that the naive idea of approximating the separating hyperplanes to enough bits of accuracy, fed into standard cutting-plane algorithms is the best one can do, at least with the bit and inner product oracles. Even in practice, algorithms for convex optimization use first-order methods that take $O(d)$ iterations to run, so they use $O(d^2)$ bits of information since each first-order query reveals at least $d$ coordinates. Our results show that this is best possible; i.e., in terms of number of bits of information used one cannot do better than what is typically done.

\subsection{Related work} Our work in this paper is also intricately related to communication complexity, since in this setting one is particularly concerned with the number and size (in bits) of messages that are sent between agents. Under the assumption that all bits of a binary representation up to a certain precision must be communicated at each step (i.e., the solver cannot ask for a specific bit of first order information without also ``paying'' for all the preceding bits), Tsitsiklis and Luo \cite{Tsitsiklis1986CommunicationCO} show a lower-bound of $d^2$ bits for communication complexity of convex optimization, and their argument turns into a lower-bound of only $d$ bits without the assumption. The authors themselves point out that the imposed restriction is `quite severe', since it effectively rules out any strategy that could make use of partial bit-information effectively. To illustrate why this might matter, consider the later stages of a convex optimization process. If one has queried many gradients, and builds a linear under-approximator of the objective function in the optimization process, one might know that in a specific region the gradient norm of the objective must be small, i.e., one has partial information about the gradient in some part of the search space already. In particular, that region of the search space is exactly where most algorithms will place the next iterate (for example, Newton's method choosing the next iterate as where the quadratic model predicts the gradient to be \textit{exactly} zero). 

In \cite{Ghadiri2024ImprovingTB}, the authors show a $\Omega(d)$ lower bound on communication complexity for the linear feasibility problem in a distributed setting, and show that their obtained lower-bounds are tight for least-squares regression and low-rank approximation problems. \cite{Vempala2019TheCC} shows a $\Omega(d^2)$ lower bound in communication complexity for exactly solving linear systems and linear programs. 

Considering these existing results, it is reasonable to conjecture that the bit-wise information complexity of convex optimization should be $\Omega(d^2)$, which we confirm in this paper.

\subsection{Future avenues}
In this paper, we obtain lower bounds on the information complexity of the feasibility problem for mixed-integer convex optimization. A natural conjecture is that similar $\Omega(d^2)$ lower bounds exist for unconstrained optimization. In fact, under the standard restriction that the objective function is $M$-Lipschitz continuous for some $M>0$, we conjecture that $$\Omega\left(2^nd^2\log\left(\frac{MR}{\eps\rho}\right)\right)$$ queries are needed to solve~\eqref{def: constrained opti} for both the bit oracle and the inner product oracle, where these queries can be made for the function subgradients, as well as function values\footnote{To make a version of the inner product oracle for function values, one can ask if the function value is bigger or smaller than a queried threshold.}, in addition to the normals of separating hyperplanes. Because of a general transfer result~\cite[Theorem 7]{Basu2023InformationCO}, and the feasibility results we establish in this paper, it suffices to consider unconstrained, pure continuous convex optimization. In particular, it suffices to establish a lower bound of $\Omega\left(d^2\log\left(\frac{MR}{\eps}\right)\right)$ on the number of queries using the bit oracle or inner product oracles, for the problem of minimizing an $M$-Lipschitz convex function over $[-R, R]^d$. In fact, it is unknown (to the best of our knowledge) if such a bound holds even for the much weaker oracle where only function values are revealed at queried points, i.e., a zeroth-order oracle. A resolution of this conjecture will provide a complete picture of the bit complexity of mixed-integer convex optimization.

\section{Proofs}

Given any $\a \in \R^k$ and $\delta \in \R$, we will use the notation $H^{\leq}(\a, \delta)$ to denote the halfspace $\{\z \in \R^k: \langle \a, \z \rangle \leq \delta\}$, and the notation $H^{=}(\a, \delta)$ to denote the hyperplane $\{\z \in \R^k: \langle \a, \z \rangle = \delta\}$. We use $\e^i \in \R^k$ to denote the $i$-th canonical vector. Finally, in the absence of integer variables (i.e., $n=0$), we simplify the notation of the family of instances of interest to $\cI_{d,R,\rho} := \cI_{0,d,R,\rho}$.

Given a set of query-response pairs of queries made by an algorithm and responses given by an oracle, we say that an instance is \textit{consistent} with these pairs if there exists a first order map under which the oracle would produce the given responses for the respective queries. In particular, this means that after the algorithm has made some number of queries and received the oracle's responses, if there are two instances with disjoint solutions that are consistent with those responses, the algorithm cannot tell these instances apart and thus cannot solve the problem. Hence, a way to prove a lower-bound on information complexity is to adversarially construct oracle responses to queries according to some strategy, and show that under this strategy after $k$ queries there will be two disjoint instances that are consistent with the observed query-response pairs; this forms the basis of our proofs. For a more detailed discussion of the validity of such adversarial arguments for constructing lower-bounds, see \cite[Appendix A]{Basu2023InformationCO}. 

\subsection{Proof of Theorem~\ref{thm:bit-mixed}}

As stated in Theorem~\ref{thm:dir-mixed}, we will prove the result for the stronger \textit{coordinate oracle}, which has permissible queries $\cQ = \{q_i: i\in \{1, ..., d\}\}$ with $q_i(\a)= \a_i$, i.e. $q_i$ is a query that returns the $i^{th}$ coordinate of the normal vector of a separating hyperplane. To prove lower bounds on the information complexity, it suffices to come up with an adversarial strategy for the oracle responses such that even after the stated number of queries, there are two disjoint instances that are both consistent with the oracle responses so far; see Section 2.1 of~\cite{Basu2023InformationCO} for a full discussion of this point.

We first design this adversarial oracle with no integer variables, i.e., $n=0$. 

For that, we label the orthants using $2^d$ strings $s_1\cdots s_d$, with $s_i\in \{-1, 1\}$; $s_1...s_d$ refers to the orthant that contains the vector $(s_1,\ldots,s_d)\in \doubleR^d$; this orthant will be denoted by $O_{s_1\cdots s_d}$. It suffices to show that even after $\frac{d^2}{16}$ queries, the oracle can respond in such a way that all the instances contained in the interior one of the orthants intersected with $[-R,R]^d$ are all consistent with these responses. Indeed, suppose the orthant indexed by $s_1\cdots s_d$ is such an orthant. 
Then the set of instances from $\cI_{d,R,\rho}$  in the interior of $[-R,R]^d \cap O_{s_1\cdots s_d}$ contain translations of all the instances in $\cI_{d,R/3, \rho}$, and one can simply repeat the argument (with appropriate translations).\footnote{In fact, translations of instances $\cI_{d,R/(2+\eps), \rho}$ for any $\eps >0$ work.}

Thus, we can iterate this argument $K$ times as long as $\rho < \frac12\cdot \frac{R}{3^K}$, i.e., we can repeat the argument at least $\Omega\big(\log\big(\frac{R}{\rho}\big)\big)$ times and still maintain $\rho < R'/2$, where $R'$ is the final radius of the $\ell_\infty$ box under consideration after $\Omega\big(\log\big(\frac{R}{\rho}\big)\big)$ iterations of the argument. Now there exist at least $2^d \geq 2$ disjoint $\ell_\infty$ balls of radius $\rho$ in $[-R', R']^d$. By the recursive invariant, these instances are consistent with all the responses of the oracle so far, and the algorithm cannot report a feasible point to both these instances since they are disjoint. Since in each recursive step, we forced the algorithm to make at least $\frac{d^2}{16}$, we obtain a lower bound of $\Omega\big(d^2\log\big(\frac{R}{\rho}\big)\big)$.

We now describe the adversarial strategy used to determine the responses of the oracle for the first $\frac{d^2}{16}$ queries so that the instances $\cI_{d,R, \rho}$ contained in the interior of one of the orthants are all consistent with these responses.

The oracle answers these queries as follows:
\begin{enumerate}
    \item Initialize a set $\mathcal{E} = \emptyset$, and a set $\cJ_{s_1...s_d} = \emptyset$ for each orthant $s_1...s_d$. 
    \item Answer $0$ for the first $\frac{d}{4}-1$ queries made in \textit{each} orthant $s_1...s_d$, no matter which queries were made. For a query made for coordinate $i$, update $\cJ_{s_1...s_d} \leftarrow \cJ_{s_1...s_d} \cup \{i\}$.\footnote{Notice that the orthants intersect at their boundary. Thus, we actually assign each point $\y \in \R^d$ to a unique orthant $O_{s_1,\ldots,s_d}$ containing it, e.g., that with lexicographically first string $s_1,\ldots,s_d$; all queries to $\y$ count towards its assign orthant and are answered accordingly.}
    \item When the $\frac{d}{4}$-th query is made in an orthant $s_1...s_d$, choose some standard basis vector $e_i \not\in \mathcal{E}$ with $i \not\in \cJ_{s_1...s_d}$.
    
    Update  $\mathcal{E} = \mathcal{E} \cup \{e_i\}$. Answer this query and all future queries in this orthant using $s_ie_i$ as the normal vector defining the separating hyperplane. 
    
\end{enumerate}

Since the number of queries is at most $\frac{d^2}{16}$, there will be at most $\frac{d}{4}$ orthants in which at least $\frac{d}{4}$ have been made, and so $|\mathcal{E}|\leq \frac{d}{4}$ above; thus, the oracle can always choose some $e_i \notin \mathcal{E}, i\notin \cJ$ in item 2. above.
Note also that $s_ie_i$ is the $i$-th coordinate vector pointing away from the origin in orthant $s_1...s_d$. 

We will show that after $\frac{d^2}{16}$ queries, there is at least one orthant such that all instances in the interior of that orthant are consistent with all query responses given by the oracle. In fact, the entire interior of the orthant is consistent with the responses of the oracle, i.e., the separating hyperplanes chosen by the oracle for its responses are all valid for the interior of the orthant (we just say that such an orthant is \emph{consistent} with the queries, dropping the word ``interior'' for brevity). The following two observations are useful: 
\begin{enumerate}
    \item If $\frac{d}{4}-1$ or fewer queries are made at points $\y^1,\ldots,\y^m$ in orthant $s_1 \cdots s_d$, the orthants not containing points in $\{\y^1, \ldots,\y^m\} + span\{\e^i: i\in \cJ_{s_1\cdots s_d}\}$ remain consistent with the answers given to these queries (recall that $\cJ_{s_1\cdots s_d}$ is the set of coordinates queried in these queries).  This is because the set of orthants not containing points in $\{\y^1, \ldots,\y^m\} + span\{\e^i: i\in \cJ_{ s_1 \cdots s_d}\}$ are precisely those labeled by strings $\hat{s}_1\cdots \hat{s}_d$ where $\hat{s}_i = - s_i$ for some $i\not\in \cJ_{s_1 \cdots s_d}$. One can therefore consider the separating hyperplane given by the normal vector $s_i \cdot \e^i$ at any of these query points $\y^1, \ldots, \y^m$ and it will be consistent with such an orthant. Hence, if $m$ queries are made in orthant $s_1 \cdots s_d$ that are all answered as $0$, at most $2^m$ orthants (those that do intersect $\{\y^1, \ldots,\y^m\} + span\{\e^i: i\in \cJ_{s_1\cdots s_d}\}$) are \textit{inconsistent} with these responses. 
    \item If $\frac{d}{4}$ or more queries are made in an orthant $s_1,\ldots,s_d$ at points $\{\y^1, \ldots,\y^m\}$, and coordinates $i\in \cJ_{s_1 \cdots s_d}$ were queried, the oracle ``commits'' to a separating hyperplane with normal vector $s_i \cdot \e^i$, for some $i\not \in \cJ_{s_1...s_d}$ for points in the orthant, i.e., in all future responses for separation oracle queries at points in this orthant, we will respond with the appropriate bit of $s_i \cdot e_i$. 
    Note that this and future responses in this orthant are consistent with the \textit{past responses} given to queries in this orthant since  $i\not \in \cJ_{s_1...s_d}$. 
\end{enumerate}

Without loss of generality, suppose that in $k$ orthants, $\frac{d}{4}$ queries each were made (as noted earlier, since the oracle commits to a separating hyperplane to use after $\frac{d}{4}$ queries in an orthant, further queries in that orthant are unhelpful), and in $\ell$ other orthants, $m_1, ..., m_\ell < \frac{d}{4}$ queries were made.

Suppose at most $\frac{d^2}{16}$ total queries are made. Hence, $k\leq \frac{d}{4}$ and $\sum_{i=1}^\ell m_i \leq \frac{d^2}{16}$ (in fact, the stronger bound $k\cdot \frac{d}{4} + \sum_{i=1}^\ell m_i \leq \frac{d^2}{16}$ holds, but the previous bounds are simpler and sufficient). Then from the queries in orthants in which more than $\frac{d}{4}$ queries were made, orthants in $H^\leq(-s_{\sigma(1)}e_{\sigma(1)}, 0), ..., H^\leq(-s_{\sigma(k)}e_{\sigma(k)}, 0)$ have been eliminated, where $\sigma$ is some map $\{1, ..., k\} \rightarrow \{1, ..., d\}$ labeling the choices of coordinate vectors $e_{\sigma(i)}$ the oracle committed to after $\frac{d}{4}$ queries were made in an orthant, 
where $\sigma$ labels the coordinates chosen by the oracle upon commitment, and
$s_{\sigma(i)}e_{\sigma(i)}$ is the committed normal in the corresponding queried orthant. Since $k\leq \frac{d}{4}$, these queries have eliminated at most a total of $\sum_{i=1}^{k} 2^{d-i} \leq 2^d - 2^{\frac{3d}{4}}$ orthants, by observation 2 above. 

For the queries made in orthants in which fewer than $\frac{d}{4}$ queries were made, if $m_i$ queries were made in such an orthant, these could have eliminated at most $2^{m_i}$ orthants by observation 1. above. Hence, these types of queries have removed at most $\sum_{i=1}^\ell 2^{m_i}$ orthants. Since we have $m_i \leq \frac{d}{4}$ and $\sum_{i=1}^\ell m_i \leq \frac{d^2}{16}$, we can consider the maximum number of orthants that could have been removed this way (no matter what the queries were) as: 
\begin{align}
\max_{m_1, ..., m_\ell}&\sum_{i=1}^\ell 2^{m_i} \nonumber \\
s.t. \phantom{xxx} m_i \leq \frac{d}{4} \phantom{x}&,\phantom{x}
\sum_{i=1}^\ell m_i \leq \frac{d^2}{16} \nonumber
\end{align}
which is maximized by setting $m_1, ..., m_{\min\{\ell,\frac{d}{4}\}} = \frac{d}{4}$, and all other $m_i = 0$
Hence, at most $\frac{d}{4} \cdot 2^{\frac{d}{4}} \leq 2^{d/2}$ orthants are eliminated by these queries. 

Therefore, in total no more than $2^d - 2^{\frac{3d}{4}} + 2^{\frac{d}{2}} \leq 2^d - 2^{\frac{d}{4}}$ orthants were eliminated. Since there are $2^d$ orthants, there are $2^{\frac{d}{4}} \geq 1$ (since $d\geq 1$) orthants that remain consistent with all answers given by the oracle. 

We now explain how to deal with integer variables, i.e., $n\geq 1$. At any query point $(\hat \x, \hat \y) \in [-R,R]^n \times [-R,R]^d$ such that $\hat \x \not\in [0,1]^n$, the oracle responds with a separating hyperplane separating $(\hat \x, \hat \y)$ from $[0,1]^n \times \R^d$. If $\hat \x \in [0,1]^n\setminus\{0,1\}^n$, the oracle says that the point $(\hat \x, \hat \y)$ is feasible, i.e., $(\hat \x, \hat \y) \in C$ (as per Definition~\ref{def: OUFOI} the oracle returns the $\0$ vector). Finally, if $\hat \x \in \{0,1\}^n$, then we ``lift'' the oracle response for $\hat \y$ from the continuous argument above to the $n+d$ dimensional space as follows. Let $\hat{\a} \in \R^d$ be the oracle response from the above continuous argument when queried with $\hat \y \in [-R,R]^d$. Let $\tilde \a \in \R^n$ be defined by $\tilde \a = \hat{\x} - \frac12 \ones$. Let $M_1 > 0$ be large enough such that with $\a = (M_1\tilde \a, \hat \a)$ we have $\langle \a, (\x,\y)\rangle < \langle \a, (\hat\x, \hat \y)\rangle$ for all $(\x,\y)$ that have been reported to be feasible by the oracle so far. Such an $M_1$ exists since there are only finitely many queried points $(\x, \y)$ (and so finitely that were reported to be feasible) and for all those points $\x$ is in the interior of $[0,1]^n$ and thus, $\langle \tilde \a, \x \rangle < \langle \tilde \a, \hat{\x} \rangle$. Geometrically, we take the halfspace in $\R^d$ that would be reported by the oracle strategy outlined for the continuous problem above and ``rotate'' it such that all the points that were reported to be feasible by the oracle are inside the rotated halfspace and thus, those feasible responses remain consistent with this halfspace. This can be done because $(\hat \x, \hat \y)$ lies on an extremal fiber of $[0,1]^n \times [-R, R]^d$. Similarly, let $M_2 > 0$ be such that such that with $\a = (M_2\tilde \a, \hat \a)$, $\langle \a, (\x, \y)\rangle < \langle \a, (\hat\x, \hat \y)\rangle$ for all $(\x, \y) \in [0,1]^n \times [-R, R]^d$ such that $\x \in \{0,1\}^n \setminus \{\hat \x\}$. Again, such an $M_2$ exists since for all such $(\x, \y)$, we have $\langle \tilde \a, \x \rangle < \langle \tilde \a, \hat{\x} \rangle$ for every $\x \in \{0,1\}^n \setminus \{\hat \x\}$ and $\|\y\|_\infty \leq R$. Set $M = \max\{M_1, M_2\}$. This ensures that the ``rotated'' halfspace contains all the fibers $\{\x\} \times[-R,R]^d $ for all $\x \in \{0,1\}^n \setminus \{\hat \x\}$. The oracle responds to the query $(\hat\x, \hat\y)$ with the separating hyperplane $\a = (M\tilde \a, \hat \a)$. The resulting halfspace satisfies 3 important conditions:
\begin{enumerate}
    \item It achieves the intended effect on the fiber $\{\hat\x\} \times [-R, R]^d$ as per the strategy for the continuous problem outlined above.
    \item It is consistent with all previously queried points that were reported to be feasible.
    \item It keeps all other fibers $\{\x\} \times[-R,R]^d $ with $\x \in \{0,1\}^n \setminus \{\hat \x\}$ intact, i.e, all those points are feasible to this halfspace.
\end{enumerate}

We now observe that if less than $2^n\cdot \frac{d^2}{16}\cdot \log_2\left(\frac{R}{2\rho}\right)$ queries are made, on at least one of the fibers $\{\x\} \times[-R,R]^d $ with $\x \in \{0,1\}^n$, less than $\frac{d^2}{16}\cdot \log_2\left(\frac{R}{2\rho}\right)$ queries have been made. By the argument for the continuous case above, there are two disjoint $\ell_\infty$ balls $B_1, B_2$ of radius $\rho$ on that fiber. Let $C_1$ be the the convex hull of all the queried points that were reported to be feasible by the oracle and $B_1$, and similarly, let $C_2$ be the the convex hull of all the queried points that were reported to be feasible by the oracle and $B_2$. Since all points $(\x,\y)$ reported to be feasible have $\x \in[0,1]^n\setminus\{0,1\}^n$, and $B_1 \cap B_2 = \emptyset$, we have $C_1 \cap C_2 \cap (\Z^n\times \R^d) = \emptyset$, i.e., $C_1$ and $C_2$ have no common mixed-integer point. All of the oracle responses are consistent with both $C_1$ and $C_2$, and thus the oracle cannot correctly report a feasible point after less than $2^n\cdot \frac{d^2}{16}\cdot \log_2\left(\frac{R}{2\rho}\right)$ queries.

\color{black}

\subsection{Proof of Theorem~\ref{thm:dir-mixed}}

As stated in Theorem~\ref{thm:bit-mixed}, we will prove the result for the stronger \textit{inner product oracle}, which has permissible queries $q_{\vv}$ for every $\vv \in \R^n \times \R^d$ with $q_\vv(\g) = \vv^T \g$, i.e., $q_\vv$ is a query that returns the inner product of the vector $\vv$ with a normal vector of a separating hyperplane. As with the proof of Theorem~\ref{thm:bit-mixed}, we first design an adversarial oracle with no integer variables, i.e., $n=0$. Moreover, we first establish a baseline $\Omega(d^2)$ lower bound, which we will improve with a $\frac{\log(R/\rho)}{\log(d)}$ factor via a recursive argument.

For the first $\lfloor d/2\rfloor$ queries, we simply report the inner product to be $0$. Notice that if the queried points are $\y^1, \ldots, \y^{\lfloor d/2\rfloor}$, and the corresponding queried directions are $\vv^1, \ldots, \vv^{\lfloor d/2\rfloor}$, there exists a hyperplane $H^{=}(\a^1,\delta_1)$ that contains $\y^1, \ldots, \y^{\lfloor d/2\rfloor}$ and whose normal vector $\a^1$ is orthogonal to $\vv^1, \ldots, \vv^{\lfloor d/2\rfloor}$. This is because we can consider the $\lfloor d/2\rfloor - 1$ difference vectors $\y^2 - \y^1, \ldots, \y^{\lfloor d/2\rfloor} - \y^1$, and there exists a nonzero vector $\a^1$ that is orthogonal to all the $2\lfloor d/2\rfloor - 1 \leq d - 1$ vectors $\y^2 - \y^1, \ldots, \y^{\lfloor d/2\rfloor} - \y^1, \vv^1, \ldots, \vv^{\lfloor d/2\rfloor}$. We normalize $\a^1$ to be of unit norm, and by switching the sign on $\a^1$, we may assume without loss of generality that $\delta_1 \geq 0$. We now consider the open (split) set $P^1 := \{\y \in \R^d: -\frac{R}{\sqrt{d}} < \langle \a^1, \y\rangle < 0\}$. Since all the queried points are contained in $H^{=}(\a^1,\delta_1)$ and the queried directions are all orthogonal to $\a^1$, the vector $\a^1$ can act as a normal vector to a separating hyperplane separating all the query points $\y^1, \ldots, \y^{\lfloor d/2\rfloor}$ from $P^1$.

In particular, all instances of closed, convex sets contained within $P^1 \cap [-R,R]^d$ are consistent with the $0$ responses of the oracle so far. 

For any subsequent queries, if they are made outside $P^1 \cap [-R,R]^d$, the oracle simply reports inner products with $\a^1$ or $-\a^1$, depending on which side of $P^1$ the queried point is, since one of these vectors can act as a separating hyperplane normal. So we focus on queried points inside $P^1 \cap [-R,R]^d$. The oracle now reports $0$ for the next $\left\lfloor \frac{d-1}{2} \right\rfloor$ such queries. By a similar argument, there is a 
$H^{=}(\a^2,\delta_2)$ that contains all these $\left\lfloor \frac{d-1}{2} \right\rfloor$ queried points, and such that the normal $\a^2$ is orthogonal to all the queried directions associated with these queried, {\em and} $\a^2$ is orthogonal to $\a^1$ (the $-1$ in $\left\lfloor \frac{d-1}{2} \right\rfloor$ guarantees we can satisfy this additional orthogonality constraint). Again we normalize $\a^2$ to be of unit norm, by possibly switching the sign on $\a^2$ we may assume without loss of generality that $\delta_2 \geq 0$, and consider the set $P^2 := \{\y \in \R^d: -\frac{R}{\sqrt{d}} < \langle \a^2, \y\rangle < 0\}$. As before, the vector $\a^2$ can act as a normal vector to a separating hyperplane separating all these latest $\left\lfloor \frac{d-1}{2} \right\rfloor$ queries from $P^2$. In particular, all instances of closed, convex sets contained within $P^2 \cap [-R,R]^d$ are consistent with the $0$ responses of the oracle for these latest queries. Therefore, the closed convex sets in $P^1 \cap P^2 \cap [-R,R]^d$ are consistent with all the answers of the oracle thus far. 

Proceeding this way, after $\frac{d^2}{8} \le \left\lfloor \frac{d}{2} \right\rfloor + \left\lfloor \frac{d-1}{2} \right\rfloor + \ldots + 1$ 
queries (assuming $d\geq 2$) we have constructed orthogonal vectors $\a^1,\a^2,\ldots,\a^k$ ($k \le d$) and their associated split sets $P^i := \{\y \in \R^d: -\frac{R}{\sqrt{d}} < \langle \a^i, \y\rangle < 0\}$ such that all convex sets in $\bar{P} := P^1 \cap P^2 \cap \ldots \cap P^k \cap [-R,R]^d$ are consistent with all the answers of the oracle thus far. We claim that $\bar{P}$ contains a copy of the smaller cube $[-\frac{R}{3d}, \frac{R}{3d}]^d$. To see that, consider the vector $\u := -\big(\frac{R}{2\sqrt{d}}\,\a^1 + \ldots + \frac{R}{2\sqrt{d}}\,\a^k\big)$, which can be thought as moving to the ``middle'' of all the split sets $P^i$. Notice that the (open) $\ell_2$-ball $B\big(\u, \frac{R}{2\sqrt{d}}\big)$ is contained in each $P^i$, since for any vector $\vv$ of length $\frac{R}{2\sqrt{d}}$, we have (using the orthogonality of the $\a^i$'s) $\ip{\a^i}{\u + \vv} = -\ip{\a^i}{\frac{R}{2\sqrt{d}} \a^i} + \ip{\a^i}{\vv} = - \frac{R}{2\sqrt{d}} +  \ip{\a^i}{\vv} \in [-\frac{R}{\sqrt{d}}, 0]$. Moreover, this ball $B\big(\u, \frac{R}{2\sqrt{d}}\big)$ is also contained in $[-R,R]^d$, since for any vector $\vv$ of length $\frac{R}{2\sqrt{d}}$ we have $\|\u + \vv\|_{\infty} \le \|\u + \vv\|_2 \le \|\u\|_2 + \|\vv\|_2 = \sqrt{d} \cdot \frac{R}{2\sqrt{d}} + \frac{R}{2\sqrt{d}} \le R$, where the equality again uses the orthogonality of the $\a^i$'s. Thus, this open ball $B\big(\u, \frac{R}{2\sqrt{d}}\big)$ is contained in $\bar{P}$, and since this ball contains the cube $Q := \u + \left[-\frac{R}{3d}, \frac{R}{3d}\right]^d$, we obtain that claim that $\bar{P}$ contains the desired cube. 

Therefore, after $\frac{d^2}{8}$ queries we were able to find a cube $Q$ where all closed convex sets in it are consistent with the queries thus far; in particular, under the assumption that $\rho < \frac{R}{6d}$, we can find two disjoint instances of $\cI_{d,R,\rho}$ in this cube consistent with all the answers, and thus fewer than $\frac{d^2}{8}$ queries are not enough to solve the feasibility problem.

From this point on, for every query not in $Q$ the oracle can simply report inner products with one of the canonical vectors $\pm \e^1,\ldots,\pm\e^d$, since one of these vectors can act as a hyperplane normal separating the queried point from (all closed convex sets in) $Q$. For queries in the cube $Q$, we can repeat this entire argument since all instances in $Q$ are simply translations of the instances in $\big[-\frac{R}{3d}, \frac{R}{3d}\big]^d$. Since each such repetition yields a cube with sides a factor $\frac{1}{3d}$ smaller than those of the previous cube, we can perform this repetition $K$ times about as long as $\rho < \frac{1}{2}\cdot \frac{R}{(3d)^K}$, i.e., for $K = \Big\lfloor \frac{\log(R/2\rho)}{\log(3d)} \Big\rfloor$ times, and each stage of the argument needs $\frac{d^2}{8}$ queries. Thus, one needs at least $\Omega\left(d^2\Big(1 + \frac{1}{\log d}\log\big(\frac{R}{\rho}\big)\Big)\right)$ queries, proving Theorem~\ref{thm:dir-mixed} in the absence of integer variables.

To handle the integer variables, i.e. $n \geq 1$, we use the same idea as in the proof of Theorem~\ref{thm:bit-mixed} where we report feasibility for query points $(\x, \y)$ with $\x \in [0,1]^n \setminus\{0,1\}^n$, and for points with $\x \in \{0,1\}^n$, we give the same responses as the oracle would for the continuous problem and then ``rotate'' the corresponding separating hyperplanes to remain consistent with all feasible points as well as other fibers.

\backmatter

\bmhead{Acknowledgments} Amitabh Basu gratefully acknowledges the support from Air Force Office of Scientific Research (AFOSR) grant FA95502510038.

\bmhead{Competing Interests}
The authors have no competing interests to declare that are relevant to the content of this article.

\bibliography{full_bib}

\end{document}